\documentclass[12pt]{amsart}
\usepackage{graphicx}
\usepackage{epstopdf}
\usepackage{url,verbatim}

\RequirePackage[colorlinks,citecolor=blue,urlcolor=blue]{hyperref}
\usepackage{breakurl}
\theoremstyle{plain}
\newcounter{mycount1}\newcounter{mycount2}\newcounter{mycount3}\newcounter{mycount}

\newenvironment{numlist}{\begin{list}{(\arabic{mycount2})}%
   {\usecounter{mycount2}\labelwidth=1cm\itemsep 0pt}}{\end{list}}

\newenvironment{Alist}{\begin{list}{\rm\MakeUppercase{\alph{mycount}}.}%
   {\usecounter{mycount}\labelwidth=1cm\itemsep 0pt}}{\end{list}}

\numberwithin{equation}{section}
\numberwithin{figure}{section}

\title[700-seat no-loss composition for the European Parliament]{A 700-seat no-loss composition\\ for the 2019 European Parliament}
\author{G. R. Grimmett}

\address{Centre for
Mathematical Sciences, Cambridge University, UK} 

\email{g.r.grimmett@statslab.cam.ac.uk}
\urladdr{\url{http://www.statslab.cam.ac.uk/~grg/}}

\author{F. Pukelsheim}
\address{Institut f\"ur Mathematik,
Universit\"at Augsburg, Germany}
\email{pukelsheim@math.uni-augsburg.de}
\urladdr{\url{https://www.math.uni-augsburg.de/emeriti/pukelsheim/}}

\author{V. Ram\'irez Gonz\'alez}
\address{Department of Applied Mathematics, University of Granada, Spain}
\email{vramirez@ugr.es}
\urladdr{\url{http://www.ugr.es/~vramirez/}}

\author{\\ W. S{\l}omczy\'nski}
\address{Institute of Mathematics,
Jagiellonian University, Krak\'ow, Poland}
\email{wojciech.slomczynski@im.uj.edu.pl}
\urladdr{\url{http://www.im.uj.edu.pl/instytut/pracownik?id=172}}

\author{K. {\.Z}yczkowski}
\address{Marian Smoluchowski Institute of Physics,
Jagiellonian University, Poland} 
\email{karol@tatry.if.uj.edu.pl}
\urladdr{\url{http://chaos.if.uj.edu.pl/~karol/}}

\begin{document}

\begin{abstract}
\emph{The following paper is 
part of the authors' response to an invitation from the  Constitutional Affairs Committee (AFCO) 
of the European Parliament to advise on mathematical methods for the  allocation of Parliamentary seats between
the 27 Member States following the planned departure of the United Kingdom in 2019. 
The authors were requested
to propose a method that respects the usual conditions of EU law,
and with the additional property that no Member State 
 (other than the UK) receives fewer that its 2014 allocation. This paper was delivered
to the AFCO on 21 August 2017, for consideration by the AFCO at its meeting 
in Strasbourg on 11 September 2017.}
\end{abstract}
\date{20 September 2017}
\maketitle
\noindent\hfil\rule{0.5\textwidth}{.4pt}\hfil

\medskip

{\leftskip=3pc\rightskip=3pc\parindent=0pt\small\parskip=6pt
{\sc Executive Summary.}
A composition for the 2019 EP is described obeying the following criteria: (1) no Member State loses any seats; (2) degressive proportionality is respected; and (3) the EP size is 700 seats. The allocation grants five base seats to every Member State, plus one seat per 890 000 citizens or part thereof, except when: (i) the no-loss criterion (1) or the degressivity criterion (2) warrants more seats, or (ii) the maximum capping imposes 96 seats. 

The approach extends to other EP sizes such as 701 seats (with allocation key 888 600), or with 710 seats (allocation key 850 000).

{\sc Abstract.}
The No-Loss Cambridge Compromise is a pragmatic method for the allocation of the seats of the 2019 European Parliament between the EU27 Member States. It has the following properties: no Member State loses any seat; it is objective and transparent, but not durable; it satisfies degressive proportionality; it is fair subject to the requirements of no-loss and degressivity; Parliament sizes as small as 700 or 710 seats may be achieved; the least possible Parliament size is 694 with current data.

}

\section{Introduction}\label{sec:int}

In January 2017, the AFCO Committee organized a workshop on the composition of the European Parliament (EP); see Grimmett et al.\ (2017). Various methods were presented for the allocation of seats between the Member States, each respecting degressive proportionality and satisfying the conditions of objectivity, fairness, durability and transparency. The Cambridge Compromise of Grimmett et al.\ (2011) may be used to achieve an EP of general size, while requiring transfers of seats between certain Member States. Other methods presented at the workshop permit each Member State to retain its 2014 seats (no-loss solutions), while requiring substantially more than 700 seats.

A more focussed question is the following. Is it possible to modify any of the proposed methods, or to devise a new method, so that the following criteria are met:
\begin{numlist}
\item 	no Member State loses any seats;
\item	degressive proportionality is respected;
\item	the size of the EP is 700 seats, or 701 or 710 seats.
\end{numlist}
This note proposes a variant of the Cambridge Compromise that meets these criteria. The variant is called the No-Loss Cambridge Compromise.

The No-Loss Cambridge Compromise operates in three stages. Stage A responds to criterion (1), Stage B to criterion (2), and Stage C to criterion (3):
\begin{Alist}
\item Allocate to each Member State its seat-count from the 2014 allocation. (For EU27, this requires a total of 678 seats.)
\item	Allocate a minimum number of further seats in order to achieve degressive propor-tionality. The ensuing seat-totals are termed minimum restrictions. (This requires a further 16 seats. There is no choice in their allocation.)
\item	Allocate the remaining 6 seats to bring the total to 700, or the remaining 16 seats to bring the total to 710, using an adapted version of the Cambridge Compromise 
\end{Alist}
Stage A is dictated by the political requirement that no Member State loses seats. Stage B is in response to the fact that the ensuing composition of the sitting EP fails to be degressively proportional. Stages A and B bind $678 + 16 = 694$ seats according to current data. Stage C has only a few additional seats at its disposal to make up the required total.

\section{Assessment according to the criteria}

The No-Loss Cambridge Compromise qualifies to be \emph{objective} and \emph{transparent}.

The issue of \emph{fairness} is more nuanced, since several Member States enjoy protected status at Stage B. The method is, however, fair provided the notion of fairness is subordinated to criteria (1) and (2) above. 

The method cannot claim to be \emph{durable}. It is an ad hoc modification devised for the incumbent parliament, and it will inevitably require modification for future parliaments.

\section{Determination of minimum restrictions (Stages A and B)}
Table 1 exhibits the 27 Member States by decreasing population figures “QMV2017”. The population figures are identical to those decreed for Council’s Qualified Majority Voting rule during the calendar year 2017; see European Council (2016). The status-quo composition is displayed in column \lq\lq 2014". Column \lq\lq RR2014" contains the representation ratios of the Member States. The representation ratio of a Member State is defined as the quotient of its population figure and its 2014 seat allocation, rounded to the nearest whole number. 

A representation ratio that is larger than a ratio of a more populous predecessor-state 
constitutes a breach of degressivity, and is marked by an asterisk (*).

For instance, France is marked because its representation ratio is larger than that of Germany (900 833 $>$ 854 838). This is a breach of degressivity: France has a smaller population than Germany, but a French deputy represents more citizens than a German deputy. A breach of degressivity might be rectified by the retraction of a suitable number of seats from the more populous state; this is, however, forbidden by criterion (1). The other solution is to augment the total of the less populous state. For France, four further \lq\lq add-on" seats are required to attain a representation ratio consistent with degressivity 
(854 636 $<$ 854 838).

By the same argument, further seats are allocated to Spain (2), The Netherlands (2), Sweden (1), Austria (1), Denmark (1) and Ireland (2). 

The adjustment of Denmark and Ireland causes collateral breaches of degressivity that are marked by a dagger ($\dagger$). They are healed by add-on seats for Finland (1), Slovakia (1) and Croatia (1). 

A total of 16 \lq\lq add-on" seats are brought to life, see Table 1. These are added to the \lq\lq 2014" seats to yield the \lq\lq 2014DP" seats. The representation ratios for the latter, as given in column \lq\lq RR2014DP", are decreasing and satisfy degressivity.

At Stage C, the degressively proportional seat allocation \lq\lq 2014DP" is passed into the No-Loss Cambridge Compromise in the form of minimum restrictions.

\section{The 700-Seat No-Loss Cambridge Compromise (Stage C)}
We describe next the No-Loss Cambridge Compromise in the current setting. The allocation proceeds in two steps. In the 2019 setting, the first step equips every Member State with 5 base seats. This utilizes $27 \times 5 = 135$ seats, and leaves $700 - 135 = 565$ seats available for the second step.

The second step invokes the ranges of seat-counts that submit to the Union’s primary law and that obey the no-loss criterion (1) and the degressivity criterion (2). Primary law demands that the seat allocation of a Member State lies between 6 and 96 seats. Criteria (1) and (2) raise the lower limit 6 to the new limit that is exhibited in the 
\lq\lq 2014DP" column in Table 1. By setting aside the five base seats from the first step, Table 2 exhibits the pertinent \lq\lq Range" that needs to be observed by the no-loss allocation.

Once the pertinent ranges have been identified, one executes the No-Loss Cambridge Compromise. For the 2019 setting it operates as follows:
\begin{quotation}
Every Member State is allocated five base seats, plus one seat per 890 000 citizens or part thereof, except when the minimum restriction warrants more seats or the maximum capping imposes fewer seats.
\end{quotation}
In numerical terms, the population figure of a Member State is divided by the divisor 890 000. The emerging quotient is rounded upwards and increased by 5, in order to account for the base seats. Let $N$ designate the whole number thus obtained. If $N$ falls into the pertinent range (see above), then $N$ signifies the seat allocation of this state (France, Italy, Spain, and Slovenia till Malta). If $N$ lies to the right of the range then the allocation becomes 96 seats (Germany). If $N$ lies to the left of the range then the allocation is given by the applicable minimum restriction (Poland till Lithuania). The key for this allocation rule is the divisor (890 000). The divisor is determined so as to deal out exactly 700 seats.

The No-Loss Cambridge Compromise results in the allocation labelled \lq\lq 2019" in Table 2. Of the 700 seat total, 678 are bound in stage A by the no-loss criterion (1), 16 are purposefully added in stage B due to the criterion of degressivity (2), and 6 are submitted to the proper allocation process in Stage C. These six seats are assigned to France (2), Italy (1), Spain (2) and Estonia (1). 

Column \lq\lq RR2019" confirms degressivity; the representation ratios are decreasing when passing from a more populous state to a less populous state. The representation ratio of a Member State is the quotient of its population figure and its 2019 seat allocation, rounded to the nearest whole number.

One may contemplate the use of a different method in Stage C, such as one of the procedures that were presented during the AFCO workshop or variants thereof. However, in 2019 only six seats are available at stage C, whence we see no compelling reason to employ a method more sophisticated than the No-Loss Cambridge Compromise. Allocations for future legislative periods of the EP will require a renewed analysis. \emph{We emphasize the simplicity of the current proposal, and the clarity of its three Stages in response to the three criteria.}

\section{Versions with 701 or 710 seats}
Art.\ 14(2) of the Treaty of Lisbon (European Union, 2012), in its peculiar wording, limits the size of the EP to 751 seats by saying that it shall not exceed seven hundred and fifty in number, plus the President. When mimicking this format, a size of 700 seats might be paraphrased as six hundred and ninety-nine in number, plus the President.

In a similar vein one may choose an EP of size seven hundred in number, plus the President. Then the total EP size would be an odd number: 701. The 701st seat would be allocated to France, with allocation key lowered to 888 600.

The No-Loss Cambridge Compromise can be applied to any parliament size of 694 or more, with current data. For a parliament of size 710 the allocation key would turn out to be 850 000 (no table included). The ten additional seats beyond 700 would increase the allocations of France (4), Italy (4) and Spain (2). This solution, also, satisfies degressivity.

However, a parliament of size much beyond 700 seats sends a different political message to the public and the media. A particular figure such as 710 would call for justification, and may underline a perception that the prime aim of the sitting parliament is to defend their incumbent mandates. In contrast, a 700 seat total is easily remembered, and conveys a more powerful message of efficiency and self-restraint.

\section{The role of rounding in degressive proportionality}
The No-Loss Cambridge Compromise has a non-standard relationship with the specification of degressive proportionality in Art.\ 1 of the European Council (2013) decision:
\begin{quotation}
The principle of degressive proportionality shall require decreasing representation ratios when passing from a more populous Member State to a less populous Member State, where the representation ratio of a Member State is defined to be the ratio of its population figure relative to its number of seats before rounding.
\end{quotation}
First of all this definition cannot be applied literally to the 2014 composition in Table 1 since no rounding is taking place.

Nor can the definition be applied literally to Table 2. The reason is that for 18 Member States (Germany, and Poland till Lithuania) the seat numbers are achieved by restriction to given ranges, and they involve no rounding step. The other nine Member States allow a distinction between before rounding and after rounding. However, the version before rounding does not interact seamlessly with the 18 Member States just mentioned: Slovenia's index 2 064 188/7.3 = 282 765 would surpass Lithuania's 262 596, thus indicating a breach of degressivity. 

Therefore the version after rounding is applied throughout. It transpires well in the 2019 setting with current data.

\vfill\eject

\begin{table}[t]
\centerline{\includegraphics[height=0.6\textheight]{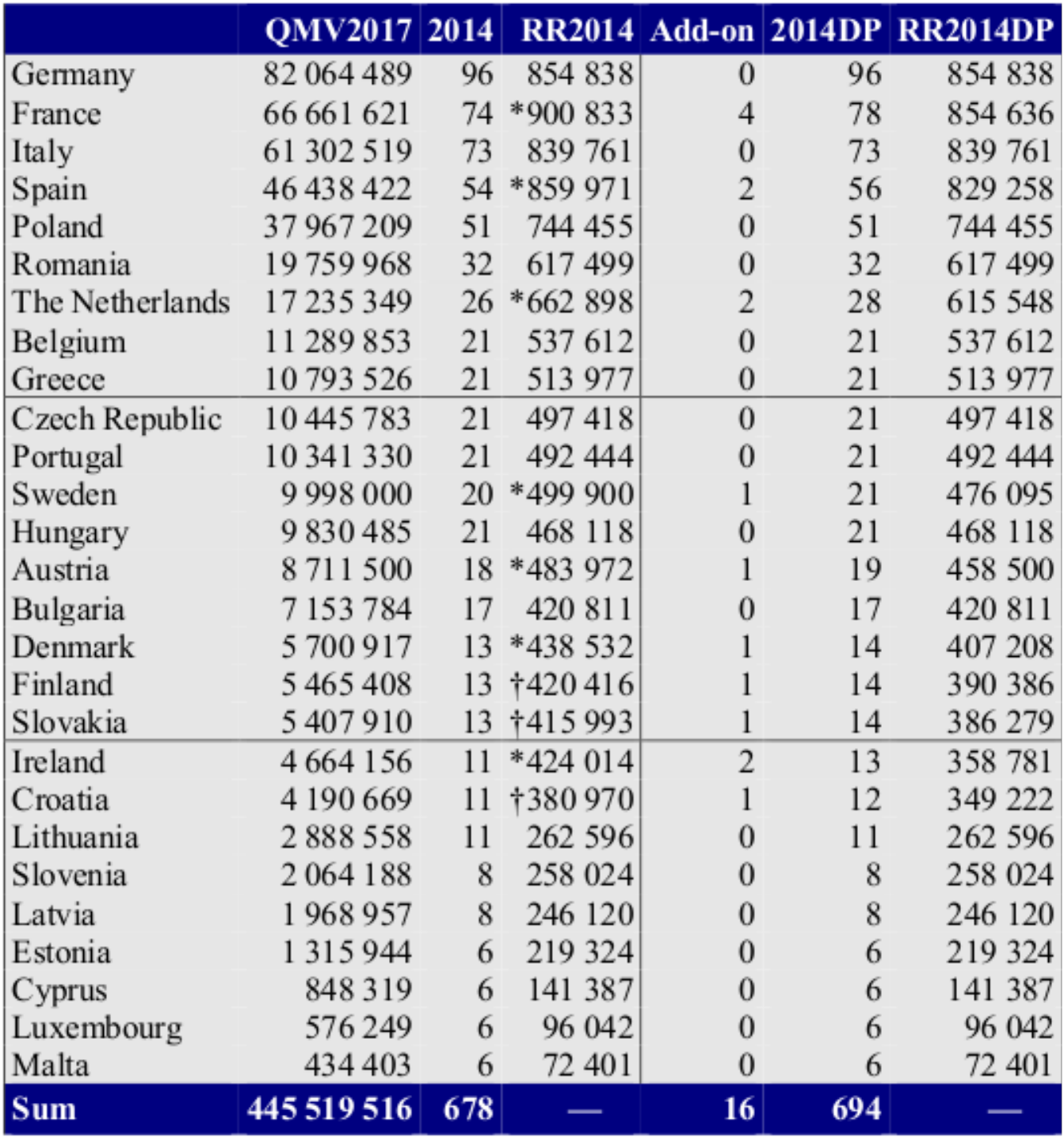}}

\caption{Augmentation of the 2014 EP composition to achieve degressivity}
\end{table}

{\footnotesize
\par\noindent{\bf Notes:}

\noindent
Representation ratio \lq\lq RR2014" of status quo \lq\lq 2014" seats:
\begin{quotation}
The representation ratios \lq\lq RR2014" are the quotient of a Member State's \lq\lq QMV2017" population figure and its \lq\lq 2014" seat allocation, rounded to the nearest whole number. A representation ratio that is larger than one of its more populous predecessors constitutes a breach of degressivity, and is marked by an asterisk (*) or, when implied by previous corrections, by a dagger ($\dagger$).
\end{quotation}

\noindent
Representation ratio \lq\lq RR2014DP" of \lq\lq 2014DP" seats:
\begin{quotation}
The \lq\lq Add-on" seats are added to the \lq\lq 2014" seats in order for the resulting “\lq\lq 2014DP" seats to achieve degressivity. The representation ratios \lq\lq RR2014DP" are the quotient of the \lq\lq QMV2017" population figures and the \lq\lq 2014DP" seats, rounded to the nearest whole number. These representation ratios are decreasing when passing from a more populous Member State to a less populous Member State.
\end{quotation}
}
\vfill\eject

\begin{table}[t]
\centerline{\includegraphics[height=0.6\textheight]{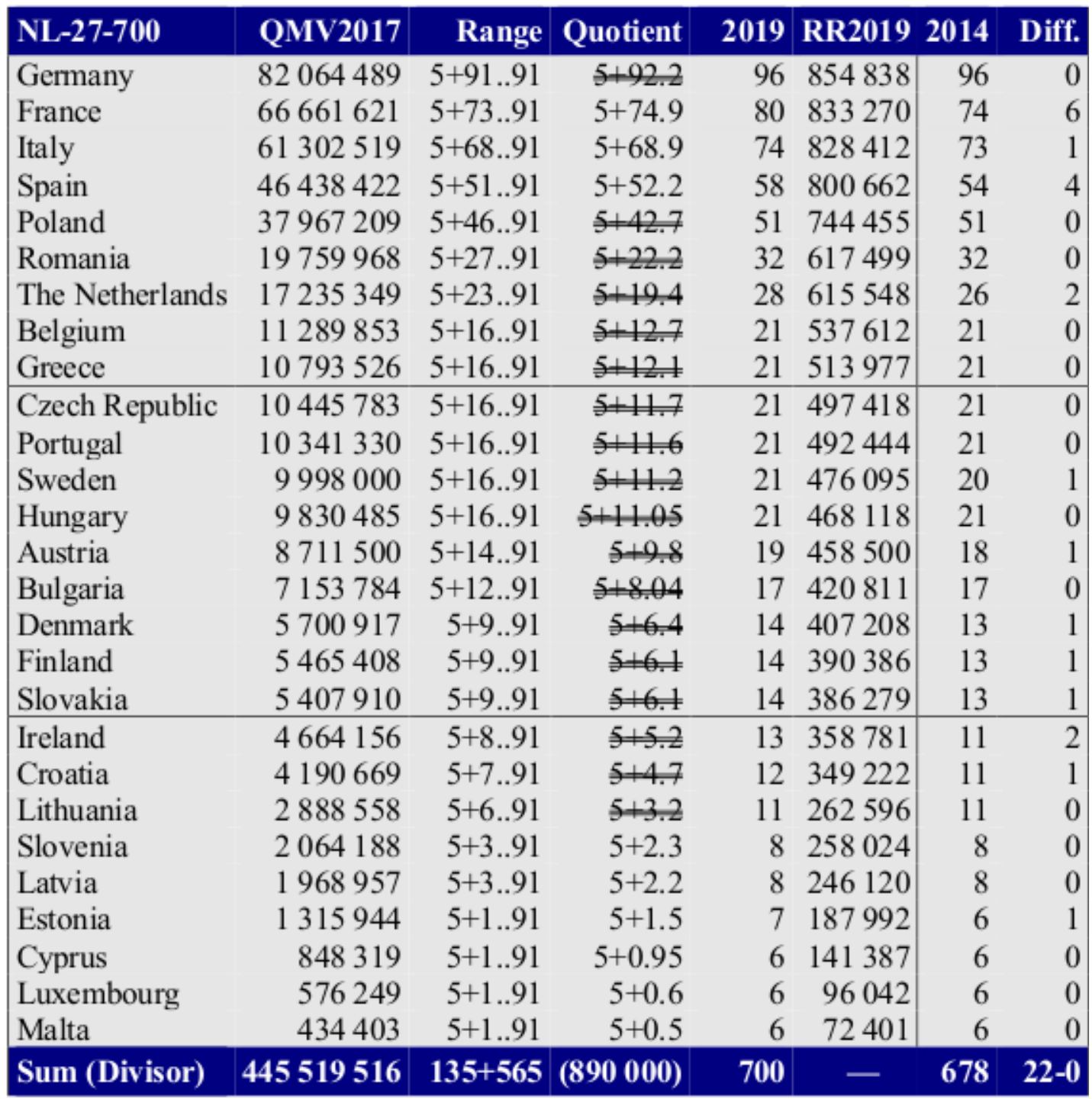}}

\caption{No-Loss Cambridge Compromise for a 2019 EP with 700 seats}\vskip-0.5cm
\end{table}

{\footnotesize
\par\noindent{\bf Notes:}

\noindent
No-Loss Cambridge Compromise \lq\lq 2019":
\begin{quotation}
Every Member State is allocated 5 base seats, plus one seat per 890,000 citizens or part thereof, except when the minimum restriction warrants more seats (Poland till Lithuania) or the maximum capping imposes fewer seats (Germany).
\end{quotation}

\noindent
Range specification:

\begin{quotation}
The \lq\lq Range" of feasible seat-counts is composed of 5 base seats plus at least a minimum number of seats so as to reach the \lq\lq 2014DP" seats in Table 1, and at most a maximum number of 91 seats so as to respect the capping at 96 seats.
\end{quotation}

\noindent
Allocation key:

\begin{quotation}
The divisor (890 000) is determined so that exactly 700 seats are dealt out.
\end{quotation}

\noindent
Verification of degressive proportionality:

\begin{quotation}
The representation ratios \lq\lq RR2019" are the quotients of the \lq\lq QMV2017" population figures and the \lq\lq 2019" seat allocations, rounded to the nearest whole number. These representation ratios are decreasing when passing from a more populous Member State to a less populous Member State.
\end{quotation}

\noindent
Column \lq\lq Diff." exhibits deviations of proposed \lq\lq 2019" seats from status quo \lq\lq 2014" seats. 
}

\section*{References}
{\small

{\hangafter=1\hangindent=3pc\parindent=0pt\parskip=6pt
European Council (2013): Decision of 28 June 2013 establishing the composition of the European Parliament (2013/312/EU). Official Journal of the European Union L 181, 29.6.2013, pp.\ 57--58 
(\url{www.uni-augs-burg.de/bazi/OJ/2013L181p57.pdf}).

\hangafter=1\hangindent=3pc
European Council (2016): Decision (EU, Euratom) 2016/2353 of 8 December 2016 amending the Council's rules of procedure. Official Journal of the European Union L 348, 12.12.2016, pp.\ 27--29 
(\url{www.uni-augsburg.de/bazi/OJ/2016L348p27.pdf}).

\hangafter=1\hangindent=3pc
European Union (2012): Consolidated version of the Treaty on European Union. Official Journal of the European Union C 326, 26.10.2012, pp.\ 13--45 (\url{www.uni-augsburg.de/bazi/OJ/2012C326p13.pdf}).

\hangafter=1\hangindent=3pc
G.\ R.\ Grimmett, J.-F.\ Laslier, F.\ Pukelsheim, V.\ Ram\'irez Gonz\'alez, R.\ Rose, W.\ S{\l}omczy\'nski, M.\ Zachariasen, K.\ {\.Z}yczkowski (2011): The Allocation Between the EU Member States of the Seats in the European Parliament -- Cambridge Compromise. Note. European Parliament, Directorate-General for Internal Policies, Policy Department C: Citizen's Rights and Constitutional Affairs, PE 432.760, March 2011 (\url{www.uni-augsburg.de/pukelsheim/2011f.pdf}).

\hangafter=1\hangindent=3pc
G.\ R.\ Grimmett, F.\ Pukelsheim, V.\ Ram\'irez Gonz\'alez, W.\ S{\l}omczy\'nski, K.\ {\.Z}yczkowski (2017): The Composition of the European Parliament. Workshop 30 January 2017. European Parliament, Directorate-General for Internal Policies, Policy Department C: Citizen's Rights and Constitutional Affairs, PE 583.117, February 2017 
(\url{www.uni-augsburg.de/pukelsheim/2017c.pdf}).

}}

\end{document}